\documentclass[12pt]{article}

\usepackage{amsmath}
\usepackage{graphicx}

\usepackage[plain]{fullpage}
\usepackage{amsfonts}
\usepackage{lastpage}
\usepackage{fancyhdr}

\begin{document}

\setlength\topmargin{-0.5in}
\setlength\headsep{0.25cm}
\setlength\textheight{9in}
\setlength\textwidth{6.5in}
\setlength\oddsidemargin{-0.25in}
\setlength\evensidemargin{-0.25in}
\setlength\parindent{0.5in}
\setlength\parskip{0.0in}

\setlength\headheight{0.5in}

\newcommand{\consVec}{{\mathbf{U}}} 
\newcommand{\primVec}{{\mathbf{V}}} 
\newcommand{\charVec}{{\mathbf{W}}} 

\newcommand{\AMat}{{\mathbf{A}}} 
\newcommand{\LMat}{{\mathbf{L}}} 

\newcommand{\dx}{{\Delta x}}
\newcommand{\dt}{{\Delta t}}
\newcommand{\dtoverdx}{{\frac{\Delta t}{\Delta x}}}

\newcommand{\fH}{{\mathbf{F}}} 
\newcommand{\fD}{{\mathbf{D}}} 
\newcommand{\fS}{{\mathbf{S}}} 
\newcommand{\fTot}{{\mathbf{\mathcal{F}}}} 

\newcommand{\stressTensor}{{\mathbf{\mathcal{S}}}} 
\newcommand{\stressTensorComp}{{\mathcal{S}}} 
\newcommand{\heatFlux}{{\mathcal{Q}}} 
\newcommand{\sT}{{s}} 
\newcommand{\hF}{{q}} 

\newcommand{\comment}[1]{}

\title{Numerical Methods for the Stochastic Landau-Lifshitz Navier-Stokes Equations}

\author{
John B. Bell, Alejandro L. Garcia, Sarah A. Williams\\
\parbox{4in}{
\footnotesize
\begin {center}
Center for Computational Sciences and Engineering\\
Lawrence Berkeley National Laboratory\\
Berkeley, California, 94720, USA
\end{center}
}}

\date{}

\maketitle

\begin{abstract}
The Landau-Lifshitz Navier-Stokes (LLNS) equations incorporate
thermal fluctuations into macroscopic hydrodynamics by using
stochastic fluxes. This paper examines explicit Eulerian
discretizations of the full LLNS equations. Several CFD approaches
are considered (including MacCormack's two-step Lax-Wendroff scheme
and the Piecewise Parabolic Method) and are found to give good
results (about $10\%$ error) for the variances of momentum and energy
fluctuations. However, neither of these schemes accurately
reproduces the density fluctuations. We introduce a conservative
centered scheme with a third-order Runge-Kutta temporal integrator
that does accurately produce density fluctuations. A variety of
numerical tests, including the random walk of a standing shock wave,
are considered and results from the stochastic LLNS PDE solver are
compared with theory, when available, and with molecular simulations
using a Direct Simulation Monte Carlo (DSMC) algorithm.
\end{abstract}

\newpage

\section{Introduction}\label{IntroSection}

Thermal fluctuations have long been a central topic of statistical
mechanics, dating back to the light scattering predictions of
Rayleigh (i.e., why the sky is blue) and the theory of Brownian
motion of Einstein and Smoluchowski~\cite{Pathria:96}. More
recently, the study of fluctuations is an important topic in fluid
mechanics due to the current interest in nanoscale flows, with
applications ranging from micro-engineering
\cite{Karniadakis:05,Ho:98,Hak:99} to molecular
biology~\cite{Alberts:02,Astumian:02,Oster:02}.

Microscopic fluctuations constantly drive a fluid from its mean
state, making it possible to probe the transport properties by
fluctuation-dissipation. This is the basis for light scattering in
physical experiments and Green-Kubo analysis in molecular
simulations. Fluctuations are dynamically important for fluids
undergoing phase transitions, nucleation, hydrodynamic
instabilities, combustive ignition, etc., since the nonlinearities
can exponentially amplify the effect of the fluctuations.

In molecular biology, the importance of fluctuations can be appreciated
by noting that a typical molecular motor protein consumes ATP at a
power of roughly $10^{-16}$ watts while operating in a background of
$10^{-8}$ watts of thermal noise power, which is likened to be ``as
difficult as walking in a hurricane is for us"~\cite{Astumian:02}.
While the randomizing property of fluctuations would seem to be
unfavorable for the self-organization of living organisms, Nature
has found a way to exploit these fluctuations at the molecular
level. The second law of thermodynamics does not allow motor
proteins to extract work from equilibrium fluctuations, yet the
thermal noise actually assists the directed motion of the protein by
providing the mechanism for overcoming potential barriers.

Following Nature's example, there is interest in the fabrication of
nano-scale devices powered by \cite{Soong:00} or constructed using
\cite{Tsong:02} so-called ``Brownian motors.'' Another application
is in micro-total-analytical systems ($\mu$TAS) or ``lab-on-a-chip''
systems that promise single-molecule detection and
analysis~\cite{Craighead:00}. Specifically, the Brownian ratchet
mechanism has been demonstrated to be useful for biomolecular
separation~\cite{Oudenaarden:99,Bader:99} and simple mechanisms for
creating heat engines driven by non-equilibrium fluctuations have
been proposed~\cite{Broeck:04, Meurs:04}. Finally, exothermic
reactions, such as in combustion and explosive detonation, can
depend strongly on the nature of thermal
fluctuations~\cite{Nowakowski:03,Lemarchand:04}.

To incorporate thermal fluctuations into macroscopic hydrodynamics,
Landau and Lifshitz introduced an extended form of the Navier-Stokes
equations by adding stochastic flux terms~\cite{Landau:Fluid}. The
Landau-Lifshitz Navier-Stokes (LLNS) equations may be written as
\begin{equation}
\consVec_t + \nabla\cdot\fH = \nabla\cdot\fD + \nabla\cdot\fS
\label{BasicLLNS_Eqn}
\end{equation}
where
\begin{equation}
\consVec = \left( \begin{array}{c} \rho \\ \mathbf{J} \\ E \\ \end{array}\right)
\end{equation}
is the vector of conserved quantities (density of mass, momentum and energy).
The hyperbolic flux is given by
\begin{equation}
\fH = \left( \begin{array}{c } \rho\mathbf{v} \\
                               \rho \mathbf{v} \cdot \mathbf{v} + P\mathbf{I} \\
                               \mathbf{v}E + P\mathbf{v} \\
                    \end{array}\right)
\end{equation}
and the diffusive flux is given by
\begin{equation}
\fD  = \left( \begin{array}{c }
                            0 \\
                            \tau \\
                            \tau \cdot \mathbf{v} + \kappa \nabla T\\
                    \end{array}\right),
\end{equation}
where $\mathbf{v}$ is the fluid velocity, $P$ is the pressure, $T$ is the
temperature, and
$\tau = \eta (\nabla \mathbf{v} + \nabla \mathbf{v}^T -\frac{2}{3} \mathbf{I} \nabla \cdot \mathbf{v})$
is the stress tensor.
Here
$\eta$ and $\kappa$ are coefficients of viscosity and thermal
conductivity, respectively, where we have assumed the bulk viscosity is zero.

The mass flux is microscopically exact but the other two flux
components are not; for example, at molecular scales heat may
spontaneously flow from cold to hot, in violation of the macroscopic
Fourier law. To account for such spontaneous fluctuations, the LLNS
equations include a stochastic flux
\begin{equation}
\fS = \left( \begin{array}{c }
                            0 \\
                            \stressTensor \\
                            \heatFlux + \mathbf{v}\cdot\stressTensor\\
                    \end{array}\right),
\end{equation}
where the stochastic stress tensor $\stressTensor$ and heat flux
$\heatFlux$ have zero mean and covariances given by
\begin{equation}
\langle \stressTensorComp_{ij}(\mathbf{r},t)\stressTensorComp_{k\ell}(\mathbf{r}',t') \rangle =
2 k_B \eta T \left( \delta^K_{ik} \delta^K_{j\ell} +
\delta^K_{i\ell}\delta^K_{jk} - {\textstyle \frac{2}{3}}
\delta^K_{ij} \delta^K_{k\ell} \right)
\delta(\mathbf{r}-\mathbf{r}') \delta(t-t'),
\end{equation}
\begin{equation}
\langle \heatFlux_i(\mathbf{r},t) \heatFlux_j(\mathbf{r}',t') \rangle = 2 k_B \kappa
T^2 \delta^K_{ij} \delta(\mathbf{r}-\mathbf{r}') \delta(t-t'),
\end{equation}
and
\begin{equation}
\langle \stressTensorComp_{ij}(\mathbf{r},t) \heatFlux_k(\mathbf{r}',t') \rangle = 0,
\end{equation}
where $k_B$ is Boltzmann's constant.
The LLNS equations have been derived by a variety of approaches (see
\cite{Landau:Fluid,Bixon:69,Fox:70,Kelly:71}) and have even been
extended to relativistic hydrodynamics~\cite{Calzetta:98}. While
they were originally developed for equilibrium fluctuations (see
Appendix A), specifically the Rayleigh and Brillouin spectral lines
in light scattering, the validity of the LLNS equations for
non-equilibrium systems has been derived \cite{Espanol:98} and
verified in molecular simulations \cite{Mansour:87,Mareschal:92}.

In this paper we investigate a variety of numerical schemes for
solving the LLNS equations.  For simplicity, we restrict our attention
to one-dimensional systems, so (\ref{BasicLLNS_Eqn}) simplifies to
\begin{equation}
\frac{\partial}{\partial t} \left( \begin{array}{c}
\rho \\ J \\ E \\ \end{array}\right) =
- \frac{\partial}{\partial x} \left( \begin{array}{c}
\rho u \\
\rho u^2 + P \\
(E+P)u \\
\end{array}\right)
+ \frac{\partial}{\partial x} \left( \begin{array}{c}
0 \\
\frac{4}{3}\eta \partial_x u\\
\frac{4}{3} \eta u \partial_x u - \kappa \partial_x T\\
\end{array}\right)
+ \frac{\partial}{\partial x} \left( \begin{array}{c}
0 \\
\sT \\
\hF + u \sT \\
\end{array}\right)
\label{OneDimLLNS_Eqn}
\end{equation}
where
\begin{eqnarray}
\langle \sT(x,t)\sT(x',t') \rangle &=& \frac{1}{\sigma^2} \int dy \int dy' \int dz \int dz'
\langle \stressTensorComp_{xx}(\mathbf{r},t)\stressTensorComp_{xx}(\mathbf{r}',t') \rangle \nonumber \\
&=& \frac{8 k_B \eta T}{3\sigma} \delta(x-x') \delta(t-t')
\end{eqnarray}
and
\begin{eqnarray}
\langle \hF(x,t) \hF(x',t') \rangle &=& \frac{1}{\sigma^2} \int dy \int dy' \int dz \int dz'
\langle \heatFlux_x(\mathbf{r},t) \heatFlux_x(\mathbf{r}',t') \rangle \nonumber\\
&=& \frac{2 k_B \kappa T^2}{\sigma} \delta({x}-{x}') \delta(t-t')
\end{eqnarray}
with $\sigma$ being the surface area of the system in the $yz$-plane.

Furthermore, we take the fluid to be a dilute gas with equation of
state $P = \rho R T$
and energy density $E = c_v \rho T + \frac{1}{2}\rho u^2$.
The transport coefficients are only functions of temperature;
for example, for a hard sphere gas $\eta = \eta_0 \sqrt{T}$ and $\kappa =
\kappa_0 \sqrt{T}$, where $\eta_0$ and $\kappa_0$ are constants. The
numerical schemes developed in this paper may readily be formulated
for other fluids. Our choice is motivated by a
desire to compare with molecular simulations (see Appendix B)
of a monatomic, hard sphere gas
(for which $R = k_B/m$ and $c_v = \frac{R}{\gamma-1}$ where $m$ is
the mass of a particle and the ratio of specific heats is $\gamma = \frac{5}{2}$).

Several numerical approaches for the Landau-Lifshitz Navier-Stokes
(LLNS) equations, and related stochastic hydrodynamic equations,
have been proposed. The most successful is a stochastic
lattice-Boltzmann model developed by Ladd for simulating solid-fluid
suspensions~\cite{Ladd:93}. This approach for modeling the Brownian
motion of particles was adopted by Sharma and
Patankar~\cite{Sharma:04} using a finite difference scheme that incorporates
a stochastic momentum flux into the
\emph{incompressible} Navier-Stokes equations. By including the
stochastic stress tensor of the LLNS equations into the lubrication
equations Moseler and Landman~\cite{Moseler:00} obtain good
agreement with their molecular dynamics simulation in modeling the
breakup of nanojets.
An alternative mesoscopic
approach to computational fluid dynamics, based on a stochastic
description defined by a discrete master equation, is proposed by
Breuer and Petruccione~\cite{Breuer:93, Breuer:94}.
They
show that the structure of the resulting system recovers the fluctuations
of LLNS.

Serrano and Espa\~nol~\cite{Serrano:01} describe a finite volume
Lagrangian discretization of the continuum equations of
hydrodynamics using Voronoi tessellation. Casting their model into
the GENERIC structure~\cite{Grmela:97} allows for the introduction
of thermal fluctuations yielding a consistent discrete model for
Lagrangian fluctuating hydrodynamics. Fabritiis et
al.~\cite{Fabritiis:02, Serrano:02} derive a similar mesoscopic,
Voronoi-based algorithm using the dissipative particle dynamics
(DPD) method. The dissipative particles follow the dynamics of
extended objects subject to hydrodynamic forces, with stresses and
heat fluxes given by the LLNS equations.

In earlier work Garcia, et al.~\cite{Garcia:87} developed a simple
finite difference scheme for the linearized LLNS equations. Though
successful, that scheme was custom-designed to solve a specific
problem; it cannot be extended readily, since it relies on special
assumptions of zero net flow and constant heat flux and would be
unstable in the more general case. Related finite difference schemes
have been demonstrated for the diffusion
equation~\cite{Alexander:02}, the ``train''
model~\cite{Alexander:05}, and the stochastic Burgers'
equation~\cite{Bell:06}, specifically in the context of Adaptive
Mesh and Algorithm Refinement hybrids that couple particle and
continuum algorithms.

In the next section we develop three stochastic PDE schemes based on
standard CFD schemes for compressible flow. The schemes are tested in
a variety of scenarios in sections \ref{NumericalTestsEqSection} and
\ref{NumericalTestsNonEqSection}, measuring spatial and time
correlations at equilibrium and away from equilibrium. Results are
compared to theoretically derived values, and also to results from
DSMC particle simulations (see Appendix B). We also examine the
influence of fluctuations on shock drift, comparing results from the
LLNS solver with DSMC simulations. The concluding section summarizes
the results and discusses future work, with an emphasis on the
issues related to using the resulting methodology as the foundation
for a hybrid algorithm.

\section{Numerical Methods}\label{MethodsSection}

The goal here is to develop an Eulerian discretization of
the full LLNS equations, representing an extension of the approach discussed
in~\cite{Bell:06} to compressible flow.
We restrict consideration here to finite-volume schemes in which all of the variables
are collocated, so that the resulting method can form the basis of a hybrid
method in which a particle description (DSMC) is coupled to the LLNS discretization.
Within this class of discretizations, our aim is to recover the
correct fluctuating statistics. In this section we develop two methods based on
CFD schemes that are commonly used for the Navier-Stokes equations.
We then introduce a specialized centered scheme designed to capture
fluctuation intensities.

\subsection{MacCormack Scheme}
\label{MacCormacksubsection}
Based on the success of the simple second-order scheme
in~\cite{Garcia:87}, we first consider MacCormack's variant of
two-step Lax-Wendroff for solving fluctuating LLNS.\footnote{A
standard version of two-step Lax-Wendroff was also considered with
similar but slightly poorer results.}  The MacCormack method is
applied in the following way:
\begin{eqnarray*}
\consVec_j^{*}  &=& \consVec_j^n - \dtoverdx \left(\fH_{j}^n-\fH_{j-1}^n\right) +
                                   \dtoverdx \left(\fD_{j+1/2}^n-\fD_{j-1/2}^n\right) \\
           &~&  \;\;\;\;\;\;     + \dtoverdx\left(\fS_{j+1/2}^n-\fS_{j-1/2}^n\right)\\
\consVec_j^{**}  &=& \consVec_j^* -\dtoverdx\left(\fH_{j+1}^*-\fH_j^*\right) +
                                   \dtoverdx\left(\fD_{j+1/2}^*-\fD_{j-1/2}^*\right) \\
                & & \;\;\;\;\;\; + \dtoverdx\left(\fS_{j+1/2}^*-\fS_{j-1/2}^*\right)\\
\consVec_j^{n+1} &=& \frac{1}{2}\left(\consVec_j^n + \consVec_j^{**}\right).
\end{eqnarray*}
Here
$\fD_{j+1/2}^n$ is a simple finite difference approximation to
$\fD$.

Straightforward evaluation of $\fS$ would be
\begin{equation}
\fS_{j+1/2}=\left( \begin{array}{c}
0 \\
\sT_{j+1/2} \\
\hF_{j+1/2} + u_{j+1/2}\sT_{j+1/2} \\
\end{array}\right),
\end{equation}
but we will see that some adjustment must be made.
The approximation to the stochastic stress tensor, $\sT_{j+1/2}$, is computed as
\begin{equation}
\label{computedStochasticStressTensor}
\sT_{j+1/2}^n = \sqrt{\frac{4k_B}{3\Delta t V_c}\,
\left(\eta_{j+1} T_{j+1} + \eta_{j} T_{j}\right)}~\Re_{j+1/2}^n
\end{equation}
where $V_c$ is the volume of a cell and the
$\Re$'s are independent, Gaussian distributed random values with zero mean and
unit variance.
The approximation to the discretized stochastic heat flux, $\hF_{j+1/2}$, is evaluated as
\begin{equation}
\label{computedStochasticHeatFlux}
\hF_{j+1/2}^n = \sqrt{\frac{k_B}{\Delta t V_c}
\left(\kappa_{j+1} (T_{j+1})^2 + \kappa_{j} (T_{j})^2\right)}~\Re_{j+1/2}^n.
\end{equation}
These same stochastic flux approximations are used in all the continuum methods presented here.

The stochastic components of the flux, $\fS_{j+1/2}^{\ell }$,
are independent, identically distributed Gaussian random variables with mean zero and
variance $\sigma^2$ for $\ell = n,*$.
Substituting this into the MacCormack scheme we find that the variance in the
flux at $j+1/2$ is given by

\begin{eqnarray*}
\left\langle \delta \left(\frac{1}{2}\fS^n + \frac{1}{2}\fS^*\right)^2 \right\rangle &=&
   \left(\frac{1}{2}\right)^2 \left\langle \delta \left(\fS^n\right)^2\right\rangle +
   \left(\frac{1}{2}\right)^2 \left\langle \delta \left(\fS^*\right)^2 \right\rangle \\
&=&   \left(\frac{1}{2}\right) \left\langle \delta \left(\fS^n\right)^2 \right\rangle\\
&=&   \frac{\sigma^2}{2} .
\end{eqnarray*}

That is, the variance in the flux is reduced to half its original magnitude by the averaging used
in the two-step MacCormack algorithm.
We correct this effect by replacing $\fS_{j+1/2}$ with $\tilde{\fS}_{j+1/2}$ =
$\sqrt{2}\fS_{j+1/2}$.  The MacCormack method we use is
\begin{eqnarray*}
\consVec_j^{*}  &=& \consVec_j^n - \dtoverdx \left(\fH_{j}^n-\fH_{j-1}^n\right) +
                                   \dtoverdx \left(\fD_{j+1/2}^n-\fD_{j-1/2}^n\right) \\
           &~&  \;\;\;\;\;\;     + \dtoverdx\left(\tilde{\fS}_{j+1/2}^n-\tilde{\fS}_{j-1/2}^n\right)\\
\consVec_j^{**}  &=& \consVec_j^* -\dtoverdx\left(\fH_{j+1}^*-\fH_j^*\right) +
                                   \dtoverdx\left(\fD_{j+1/2}^*-\fD_{j-1/2}^*\right) \\
                & & \;\;\;\;\;\; + \dtoverdx\left(\tilde{\fS}_{j+1/2}^*-\tilde{\fS}_{j-1/2}^*\right)\\
\consVec_j^{n+1} &=& \frac{1}{2}\left(\consVec_j^n + \consVec_j^{**}\right).
\end{eqnarray*}
%

\subsection{Piecewise Parabolic Method}
\label{PPMsubsection}
In~\cite{Bell:06} a piecewise linear second-order Godunov scheme was shown to
be effective for solving the
fluctuating Burgers' equation.
We considered two versions of higher-order Godunov methods for the LLNS, a piecewise
linear version \cite{Colella-MUSCL} and the
Piecewise Parabolic Method (PPM)
introduced in \cite{Colella:84}.
The PPM algorithm, based on the direct Eulerian version
presented in \cite{Miller:02},
produced considerably better results than the piecewise linear scheme.
Since our goal
is to preserve fluctuations, we do not limit slopes and we do not
include discontinuity detection in the algorithm.

For this scheme the hyperbolic terms of the LLNS equations
are considered in terms of hydrodynamic and local characteristic variables.
In hydrodynamic variables we have
\begin{equation}
\frac{\partial}{\partial t}\primVec + \AMat \frac{\partial}{\partial x}\primVec = 0,
\end{equation}
where
\begin{equation}
\primVec_j =  \left( \begin{array}{c}
\rho_j \\
u_j \\
P_j
\end{array} \right).
\end{equation}

The local characteristic variables are interpolated via a
fourth-order scheme to the left ($-$) and right ($+$) edges of each
cell:
\begin{equation}
\charVec_{j,\pm}^n = \frac{7}{12}(\LMat_j\primVec_j
                    + \LMat_j\primVec_{j\pm1})
                    - \frac{1}{12}(\LMat_j\primVec_{j\mp1}
                    + \LMat_j\primVec_{j\pm2}),
\end{equation}
where $\LMat_j$ is the matrix whose rows are the left eigenvectors of $\AMat$
evaluated at $\primVec_j$.

These values, together with the cell-centered value $\charVec_j^n = \LMat_j\primVec_j$, are used
to construct a parabolic profile $\charVec_{j,k}(\theta)$ for each characteristic variable $k$ in each
cell,
\begin{equation}
\charVec(\theta) = \charVec_{j,-} + \theta\Delta \charVec_j + \theta(1-\theta)\charVec_{j6},
\end{equation}
where
\begin{eqnarray*}
\theta &=& \frac{x-(j-\frac{1}{2})\Delta x}{\Delta x},\\
\Delta \charVec_{j}^n &=& \charVec_{j,+}^n - \charVec_{j,-}^n,\textrm{ and}\\
\charVec_{j6}^n &=& 6(\charVec_j^n - \frac{1}{2}(\charVec_{j,+}^n + \charVec_{j,-}^n)).
\end{eqnarray*}

Time-centered updates are based on the sign of each local
characteristic wavespeed, $\lambda_{j,k}$:
\[
\charVec_{j,\pm,k}^{n+1/2} =\left\{\begin{array}{ll}
    \frac{1}{\nu_{j,k}} \int_{\pm\frac{1}{2}-\nu_{j,k}}^{\pm\frac{1}{2}} \charVec_{j,k}(\theta) \,d\theta, & \pm \lambda_{j,k} > 0 \\
    \charVec_{j,\pm,k}^{n} & \textrm{otherwise}
       \end{array}\right.
\]
where $\nu_{j,k} = \lambda_{j,k} \dtoverdx$.

Finally, the time-centered values are transformed back into primitive variables and
used as inputs to a Riemann problem at each cell edge.
We use the
approximate Riemann solver discussed in~\cite{Colella:85}.
This approach iterates the phase space solution in the $u-p$ plane, approximating
the rarefaction curves by the Hugoniot locus. The overall approach is able to handle
strong discontinuities and is second-order in wave strength.

Approximations to the viscous and stochastic flux terms are
discussed in section \ref{MacCormacksubsection}.  For our PPM algorithm we
center the viscous update in time, so that the complete update is as follows:
\begin{eqnarray}
\consVec_j^{*}   &=& \consVec_j^{n} -\dtoverdx\fH_j^n +\dtoverdx(\fD_j^n + \tilde{\fS}_j^n)\\
\consVec_j^{n+1} &=& \consVec_j^{n} -\dtoverdx\fH_j^n +
               \frac{1}{2}\left(\dtoverdx\right)\left(\fD_j^n + \tilde{\fS}_j^n + \fD_j^{*} + \tilde{\fS}_j^{*}\right).
\end{eqnarray}
As discussed in section \ref{MacCormacksubsection}, for the PPM scheme we
use the adjusted stochastic flux approximation
$\tilde{\fS}_{j}$ = $\sqrt{2}\fS_{j}$, since the averaging in the
time-centering reduces the variance in the flux to
half its original magnitude.

\subsection{Variance-preserving third-order Runge-Kutta}
\label{RK3subsection}

Equilibrium tests, presented in detail in the next section, show
that neither stochastic version of the traditional numerical methods
discussed above accurately represents the fluctuations in the LLNS
equations. The principal difficulty arises because there is no
stochastic forcing term in the mass conservation equation.
Accurately capturing density fluctuations requires that the
fluctuations be preserved in computing the mass flux. Another key
observation is that the representation of fluctuations in the above
schemes is also sensitive to the time step, with extremely small
time steps leading to somewhat improved results. This suggests that
temporal accuracy also plays a significant role in capturing
fluctuations. Based on these observations we have developed a new
discretization aimed specifically at capturing fluctuations in the
LLNS equations. The method is based on a third order Runge-Kutta
temporal integrator (RK3) combined with a centered discretization of
hyperbolic and diffusive fluxes.

The RK3 discretizaton can be written in the following three-stage form:
\begin{eqnarray}
\consVec_j^{n+1/3} &=& \consVec_j^n - \dtoverdx(\fTot_{j+1/2}^n-\fTot_{j-1/2}^n) \label{RK3a}\\
\consVec_j^{n+2/3} &=& \frac{3}{4}\consVec_j^n + \frac{1}{4}\consVec_j^{n+1/3} -
                          \frac{1}{4}\left(\dtoverdx\right)(\fTot_{j+1/2}^{n+1/3}-\fTot_{j-1/2}^{n+1/3}) \label{RK3b}\\
\consVec_j^{n+1}   &=& \frac{1}{3}\consVec_j^n + \frac{2}{3} \consVec_j^{n+2/3} -
                           \frac{2}{3}\left(\dtoverdx\right)(\fTot_{j+1/2}^{n+2/3}-\fTot_{j-1/2}^{n+2/3}), \label{RK3c}
\end{eqnarray}
where $\fTot = -\fH + \fD + \fS$.

Combining the three stages, we can write
\begin{equation*}
\consVec_j^{n+1} = \consVec_j^n - \dtoverdx \left[
                       \frac{1}{6}(\fTot_{j+1/2}^n-\fTot_{j-1/2}^n ) +
                       \frac{1}{6}(\fTot_{j+1/2}^{n+1/3}-\fTot_{j-1/2}^{n+1/3} ) \\
                       + \frac{2}{3}(\fTot_{j+1/2}^{n+2/3}-\fTot_{j-1/2}^{n+2/3} ) \right].
\end{equation*}

The stochastic components of the flux, $\fS_{j+1/2}^{n+ \ell }$
are independent, identically distributed Gaussian random variables with mean zero and
variance $\sigma^2$
for $\ell = 0,\frac{1}{3},\frac{2}{3}$.
Substituting this into the combined update we find that the variance
in the flux at $j+1/2$ is given by
\begin{eqnarray*}
&&\langle \delta (\frac{1}{6}(\fS_{j+1/2}^0) + \frac{1}{6}(\fS_{j+1/2}^{1/3}) +
\frac{2}{3}(\fS_{j+1/2}^{2/3})) ^2 \rangle\\
  &=& \left(\frac{1}{6}\right)^2 \langle (\delta \fS_{j+1/2}^0)^2 \rangle +
   \left(\frac{1}{6}\right)^2 \langle (\delta \fS_{j+1/2}^{1/3})^2 \rangle +
   \left(\frac{2}{3}\right)^2 \langle (\delta \fS_{j+1/2}^{2/3})^2 \rangle\\
  &=&
   \frac{\sigma^2}{2} .
\end{eqnarray*}

Thus, in the course of the RK3 algorithm, the variance in the flux
is reduced to half its original magnitude, so again we replace
$\fS_{j+1/2}$ by $\tilde{\fS}_{j+1/2}$ = $\sqrt{2}\fS_{j+1/2}$, as discussed
in section \ref{MacCormacksubsection}, and compute equations
(\ref{RK3a}-\ref{RK3c}) using $\fTot = -\fH + \fD + \tilde{\fS}$.

However, this treatment does not directly affect the fluctuations 
in density, since $\fS$ does not appear in the continuity equation. 
We can correct this effect via a special interpolation scheme:
by augmenting the variance to compensate for the density reduction
arising from the temporal averaging, the fluctuations are preserved in the
mass flux computation.

We interpolate $J$ (and the other conserved quantities)
from cell-centered values:
\begin{equation}
J_{j+1/2} = \alpha_1(J_j + J_{j+1}) - \alpha_2(J_{j-1} + J_{j+2}),
\end{equation}
where
\begin{eqnarray}
\alpha_1 &=& (\sqrt{7}+1)/4 \textrm{ and}\\
\alpha_2 &=& (\sqrt{7}-1)/4 .
\end{eqnarray}
Then in the case of constant $J$ we have exactly $J_{j+1/2} = J$ and
$\langle \delta J_{j+1/2}^2\rangle = 2\langle \delta J^2\rangle$, as
desired; the interpolation is consistent and compensates for the
variance-reducing effect of the multi-stage Runge-Kutta algorithm.
The interpolation formula is similar to the PPM spatial construction
except in the PPM construction $\alpha_1 = 7/12$ and $\alpha_2 =
1/12$. Tests based on these alternative weights produced results
intermediate to the RK3 scheme and the PPM scheme. We also
considered interpolation of primitive variables but found that
interpolation based on primitive variables led to stable but
undamped oscillatory behavior. Finally, the diffusive terms $\fD$
are discretized with standard second-order finite difference
approximations.

\subsection{Boundary Conditions} \label{BoundaryTreatmentsSubsection}
In sections \ref{NumericalTestsEqSection} and \ref{NumericalTestsNonEqSection} we consider test problems for the various PDE algorithms on either a
periodic computational domain, a computational domain bounded by thermal
walls, or a computational domain bounded by infinite reservoirs. 
Boundary conditions are implemented using ghost cells. For the periodic and reservoir boundaries, it is straightforward to determine the ghost cell data.

For the case of thermal walls, in addition to ghost cells we also use a one-sided finite difference formulation to approximate $u_x$ and $T_x$ in the
calculation of the diffusive flux. The treatment of the hyperbolic flux at thermal walls varies by method.

For thermal wall boundaries in MacCormack, conserved quantities are reflected across the boundaries of the domain.  The temperature in the ghost cells
is determined by linear extrapolation, and the no-flow condition is enforced by setting the velocity terms of the hyperbolic flux to zero within the
ghost cells.

For thermal wall boundaries in PPM, ghost cells are populated by reflecting primitive variable values across the domain boundaries, and the temperature
in the ghost cells is determined by linear extrapolation. The PPM routine takes as input the cell-centered primitive variable data and returns a Riemann
solution at each cell edge.  On the domain boundaries, we modify these Riemann solutions by enforcing fixed wall temperature (i.e.,
the pressure at the wall is taken to be a function of the fixed wall temperature) before computing the hyperbolic flux across each edge.

For thermal wall boundaries in RK3, conserved quantities are reflected across the boundaries of the domain and then interpolated onto cell edges. At the
domain boundaries we employ a Riemann solver, which ensures that the boundary treatment respects characteristic compatibility relations at the physical
boundaries. At the physical boundaries, the primitive variable values derived from the conserved-quantity interpolants are modified to enforce zero
velocity and fixed wall temperature. This vector of primitive variables provides the input to the Riemann problem on the interior side of the boundary. The
input to the Riemann problem on the exterior side of the boundary is the reflection of the interior input data. The treatment of reservoir
boundaries is similar. However, ghost cells are populated with reservoir data, wall conditions are not enforced, and the input to the Riemann problem on the
exterior side of the boundary is the reservoir data.

\section{Numerical Tests -- Equilibrium}\label{NumericalTestsEqSection}

This section presents results from a variety of scenarios in which
the three schemes described above were tested. The physical domain
is chosen to be compatible with DSMC particle simulations; see
Table~\ref{DSMC_Table} for the system's parameters and Appendix B
for a description of DSMC. The domain is partitioned into 40 cells
of equal size $\Delta x$ and hyperbolic and diffusive stability
constraints determine the maximum time step $\Delta t$:
\begin{equation}
( |u| + c_s )\dtoverdx \leq 1,
\end{equation}
\begin{equation}
\textrm{max}\left(\frac{4}{3}\frac{\overline{\eta}}{\overline{\rho}},
\frac{\overline{\kappa}}{\overline{\rho} c_v}\right)
\frac{\Delta t}{\Delta x^2} \,
\leq \frac{1}{2},
\end{equation}
where the sound speed $c_s=\sqrt{{\gamma
\overline{P}}/{\overline{\rho}}}$,
$\overline{\eta}=\eta(\overline{T})$, and
$\overline{\kappa}=\kappa(\overline{T})$; the overline indicates
reference values (e.g., equilibrium values around which the system
fluctuates). For the reference state (Argon at STP) and a cell width
of $\Delta x \approx 10^{-6}$ cm the time step used was $\Delta t =
10^{-12}$ s.

\begin{table}
  \centering
  \begin{tabular}{|l|c||l|c|}
    \hline
    Molecular diameter (Argon) & $3.66 \times 10^{8}$ & Molecular mass (Argon) & $6.63 \times 10^{23}$ \\
    Reference mass density & $1.78\times 10^{-3}$ & Reference temperature & 273 \\
    Sound speed & 30781 & Specific heat $c_v$ & $3.12 \times 10^{6}$ \\
    System length & $1.25 \times 10^{-4}$ & Reference mean free path & $6.26  \times 10^{-6}$ \\
    System volume & $1.96 \times 10^{-16}$ & Time step & $1.0 \times 10^{-12}$\\
    Number of cells & 40 & Number of samples & $10^7$ \\
    Number of DSMC particles & 5265  & DSMC collision grid size & $3.13 \times 10^{-6}$ \\
    \hline
  \end{tabular}
  \caption{System parameters (in cgs units) for simulations of a dilute gas in a periodic domain.}\label{DSMC_Table}
\end{table}

\subsection{Variances at equilibrium}\label{SpatialVariancesSection}

The first benchmark for our numerical schemes is recovering the
correct variance of fluctuations for a system at equilibrium. For
this initial test problem, we take a periodic domain with zero net
flow and constant average density and temperature.
Similar
results, not presented here,  were obtained for the case of constant non-zero net flow. The
variances are computed in 40 spatial cells from $10^7$ samples and
then averaged over the cells.

Table \ref{VarEquilResultsTable} compares the theoretical variances
(see Appendix A) with those measured in the three stochastic PDE
schemes and the DSMC particle simulation. The MacCormack and PPM
schemes do relatively poor job ($9-16\%$ error) for the variances
of density and energy. Better PPM results are obtained by decreasing
our value of $\Delta t$ by a factor of 10, to $10^{-13}$. However,
it is not desirable to run simulations at such a small time step.
Only the third-order Runge-Kutta integrator generates the correct
variance of density and energy while advancing with time steps near
the stability limit.

\begin{center}
\begin{table}[h!]
\parbox{5in}{
\caption{Variance in conserved quantities at equilibrium (computed values are accurate to approximately 0.1\%).\label{VarEquilResultsTable}
}}
\begin{tabular}{l c | c | c}
\hline
\hline
                    & $\langle \delta \rho^2 \rangle$ & $\langle \delta J^2
                    \rangle$ & $\langle \delta E^2 \rangle$ \\
Exact value                 & $2.35 \times 10^{-8}$  & 13.01               & $2.87 \times 10^{10}$\\
MacCormack scheme & $2.01 \times 10^{-8}$  & 13.31               & $2.61 \times 10^{10}$\\
Piece-wise Parabolic Method& $1.97 \times 10^{-8}$ & 13.27               & $2.58 \times 10^{10}$\\
Runge-Kutta ($3^\mathrm{rd}$ order) & $2.32 \times 10^{-8}$ & 13.65                & $2.87 \times 10^{10}$\\
Molecular simulation (DSMC)       & $2.35 \times 10^{-8}$  & 13.21                & $2.79 \times 10^{10}$\\
\hline
Percentage difference (MacCormack)      & $-14.3\%$  & $2.3\%$ & $-9.3\%$\\
Percentage difference (PPM)& $-16.0\%$& $2.0\%$ & $-10.3\%$\\
Percentage difference (RK3)& $-1.3\%$  & $4.9\%$  & $-0.1\%$\\
Percentage difference (DSMC)      & $0.0\%  $  & $1.6\%$ & $-3.1\%$\\
\hline
\end{tabular}
\end{table}
\end{center}

\subsection{Spatial correlations at equilibrium}\label{SpatialCorrelationSection}

\begin{figure}[h]
\begin{center}
$\begin{array}{c@{\hspace{0.5in}}c}
\includegraphics[width=3.25in]{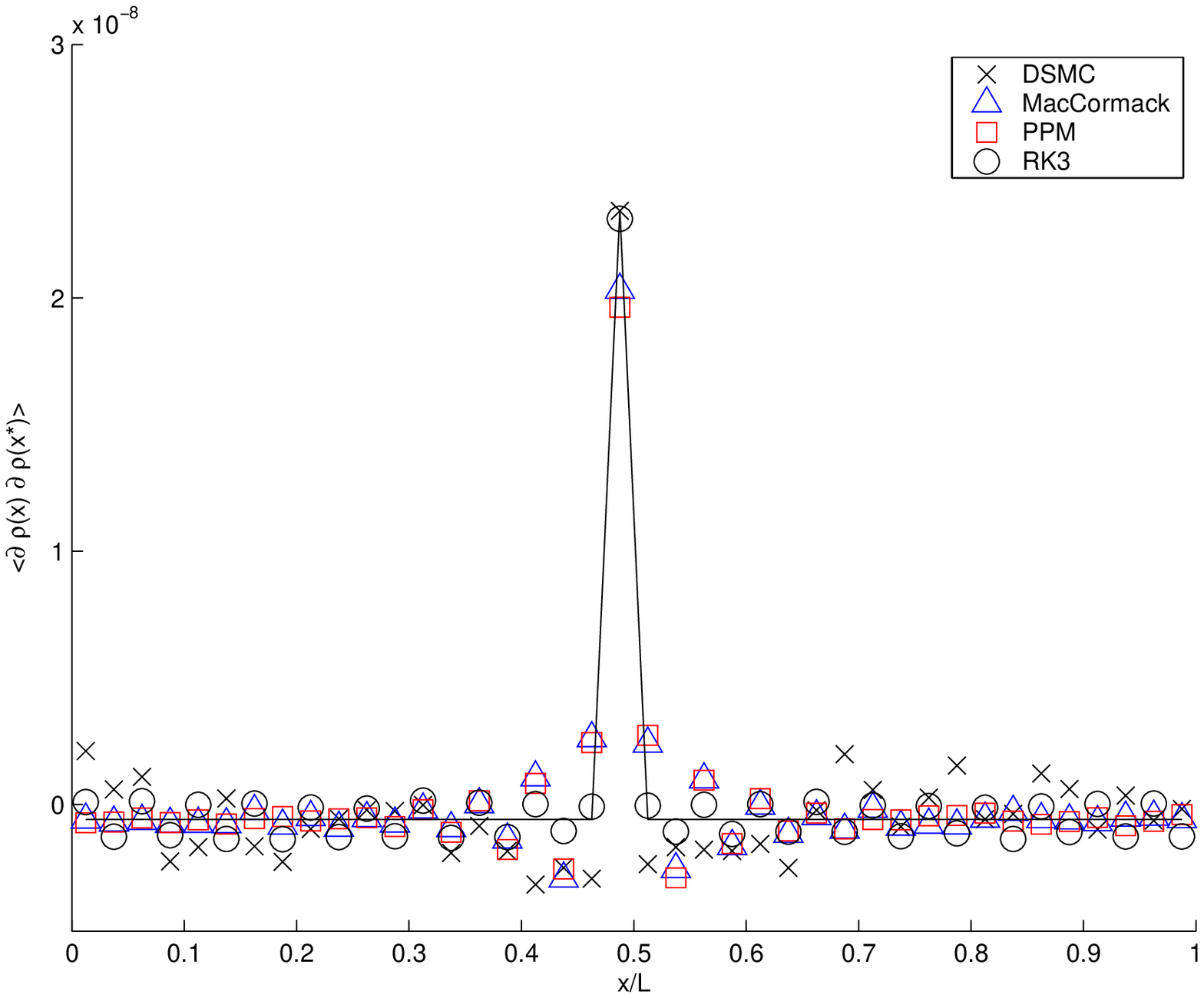}

&
\includegraphics[width=3.25in]{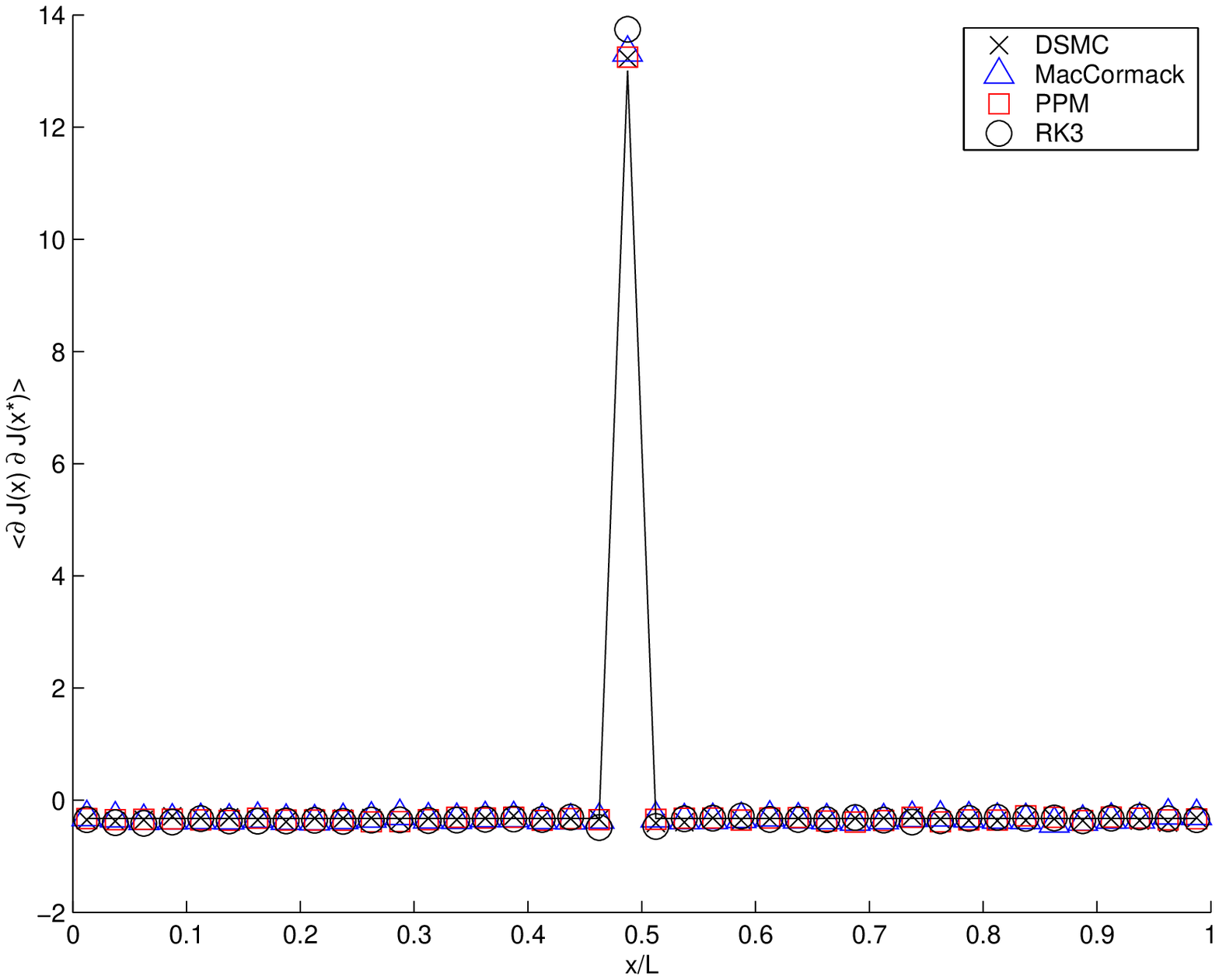} \\ [0.4cm]
\parbox{3.25in}{
\caption{\label{SpaceCorrRhoRhoFig} Spatial correlation of density
  fluctuations.  Solid line is$\langle \delta \rho_i \delta \rho_j \rangle = \langle \delta
\rho^2 \rangle \delta^K_{i,j}$ (see equations (\ref{EquilibCorrRhoFiniteEqn},
\ref{EquilibVarRhoEqn})).
}} &
\parbox{3.25in}{\caption{\label{SpaceCorrMomMomFig}
Spatial correlation of momentum fluctuations.
Solid line is$\langle \delta J_i \delta J_j \rangle = \langle \delta
J^2 \rangle \delta^K_{i,j}$ (see equations (\ref{EquilibVarMomEqn},
\ref{EquilibCorrMomFiniteEqn})). }}
\end{array}$
\end{center}
\label{SpaceCorrFigA}
\end{figure}

\begin{figure}[h]
\begin{center}
$\begin{array}{c@{\hspace{0.5in}}c}
\includegraphics[width=3.25in]{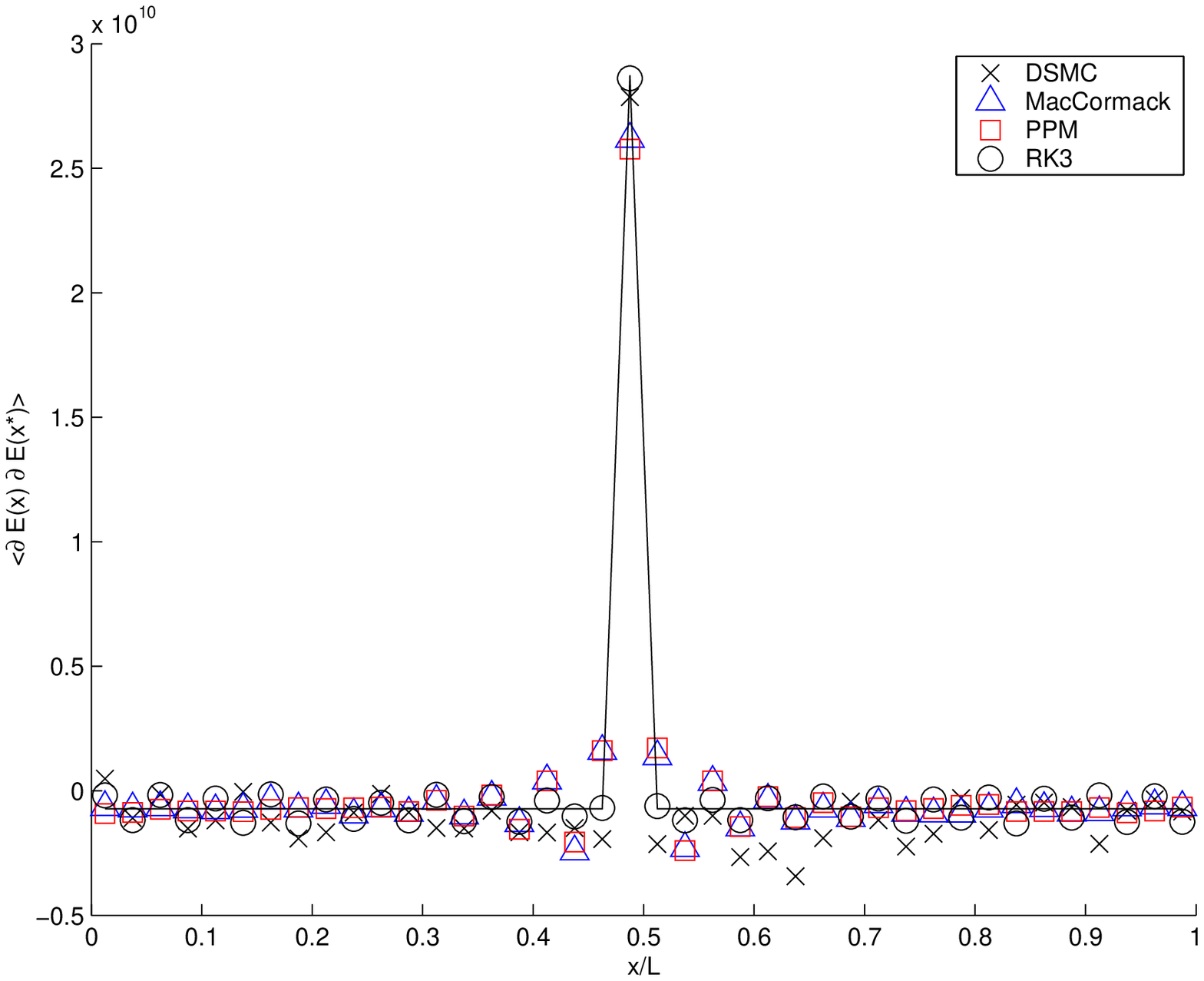} &
\includegraphics[width=3.25in]{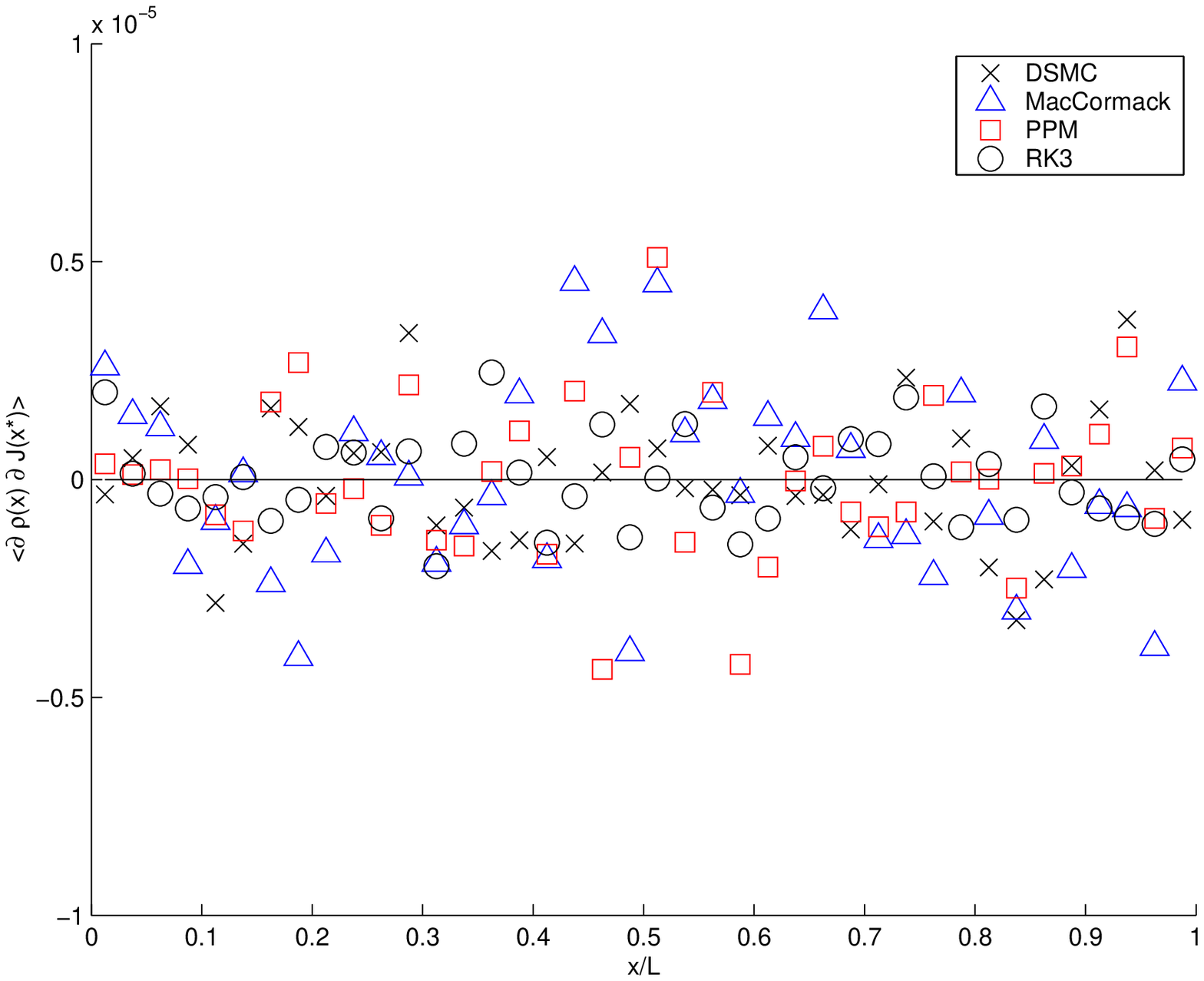} \\ [0.4cm]
\parbox{3.25in}{\caption{
\label{SpaceCorrEnrEnrFig}
Spatial correlation of energy fluctuations.
Solid line is$\langle \delta E_i \delta E_j \rangle = \langle \delta
E^2 \rangle \delta^K_{i,j}$ (see equations (\ref{EquilibVarEnrEqn},
\ref{EquilibCorrEnrFiniteEqn})).
}} &
\parbox{3.25in}{\caption{
\label{SpaceCorrRhoMomFig}
Spatial correlation of density-momentum fluctuations.}}
\end{array}$
\end{center}
\label{SpaceCorrFigB}
\end{figure}

Figures \ref{SpaceCorrRhoRhoFig}--\ref{SpaceCorrEnrEnrFig}
depict the spatial correlation of conserved variables,
that is, $\langle \delta \rho_j \delta \rho_{j^*} \rangle$,
$\langle \delta J_j \delta J_{j^*} \rangle$, and
$\langle \delta E_j \delta E_{j^*} \rangle$,
where $j^*$ is located at the center of the domain.
These figures show results computed by the MacCormack, PPM, and RK3
schemes, along with the theoretical values of the correlations (see
Appendix A) and molecular simulation data (see Appendix B). For the
MacCormack and PPM schemes the spatial correlations of density
fluctuations and energy fluctuations have significant spurious
oscillations near the correlation point (see Figs.~\ref{SpaceCorrRhoRhoFig} and
\ref{SpaceCorrEnrEnrFig}). All three schemes do well in reproducing
the expected correlations of momentum fluctuations. Figure
\ref{SpaceCorrRhoMomFig} depicts $\langle \delta \rho_j \delta
J_{j^*} \rangle$, which has a theoretical value of zero since the
net flow is zero; all three schemes correctly reproduce this result.

\subsection{Time correlations at equilibrium}\label{TimeCorrelationSection}

The time correlation of density fluctuations is of interest
because its temporal Fourier transform
gives the spectral density, which is measured experimentally from
light scattering spectra~\cite{Berne:00,Boon:91}.
From the LLNS equations, this time correlation can be written as
\begin{eqnarray}
\frac{\langle \delta \rho(w,t) \delta \rho(w,t+\tau) \rangle}
{\langle \delta \rho^2(w,t) \rangle} &=& \left( 1 - \frac{1}{\gamma}
\right)\exp\{ -w^2 D_T \tau \}
+ \frac{1}{\gamma} \exp\{ -w^2 \Gamma \tau \} \cos( c_s w \tau ) \nonumber\\
&~&\qquad + \frac{3\Gamma - D_v}{\gamma^2 c_s} w \exp\{ -w^2 \Gamma
\tau \} \sin( c_s w \tau )
\label{TimeCorrelationEqn}
\end{eqnarray}
where $w = 2\pi n /L$ is the wavenumber, $\gamma = c_p/c_v$ is the ratio of
specific heats, $D_T = \kappa/\overline{\rho}c_v$ is the thermal
diffusivity, $D_v = \frac{4}{3}\eta/\overline{\rho}$ is the
longitudinal kinematic viscosity, $c_s$ is the sound speed, and
$\Gamma = \frac{1}{2}[D_v + (\gamma-1)D_T]$ is the sound attenuation
coefficient.

In our numerical calculations the density is represented by cell
averages $\rho_i, i~=~1,\ldots,M_c$, and the time correlation is
estimated from the mean of $N$ samples,
\begin{equation}
\langle \delta \rho(w,t) \delta \rho(w,t+\tau) \rangle_N =
\frac{1}{N}\sum_{\mathrm{samples}}^N R(t) R(t+\tau)
\end{equation}
with
\begin{equation}
R(t) = \frac{1}{M_c}\sum_{i=1}^{M_c} \rho_i \sin( 2 \pi n x_i / L ).
\end{equation}
We have
\begin{equation}
\langle \delta \rho(w,t) \delta \rho(w,t+\tau) \rangle
=\lim_{N\rightarrow\infty} \langle \delta \rho(w,t) \delta \rho(w,t+\tau)
\rangle_N.
\end{equation}
From the above we find the normalization of the theoretical result
may be expressed as
\begin{eqnarray}
\langle \delta \rho^2(w,t) \rangle = \langle R(t)^2 \rangle &=&
\frac{1}{M_c^2}\sum_{i=1}^{M_c}\sum_{j=1}^{M_c} \langle \delta
\rho_i \delta \rho_j \rangle
\sin( 2 \pi n x_i / L ) \sin( 2 \pi n x_j / L )\nonumber\\
&=& \frac{\langle \delta \rho^2 \rangle}{2 M_c}.
\end{eqnarray}
We restrict our attention to the lowest wavenumber (i.e., $n=1$)
because for the system sizes we consider the theoretical result,
(\ref{TimeCorrelationEqn}), is not accurate at short wavelengths due
to mean-free-path corrections.

\begin{figure} [h]
\centering\includegraphics[width=5in]{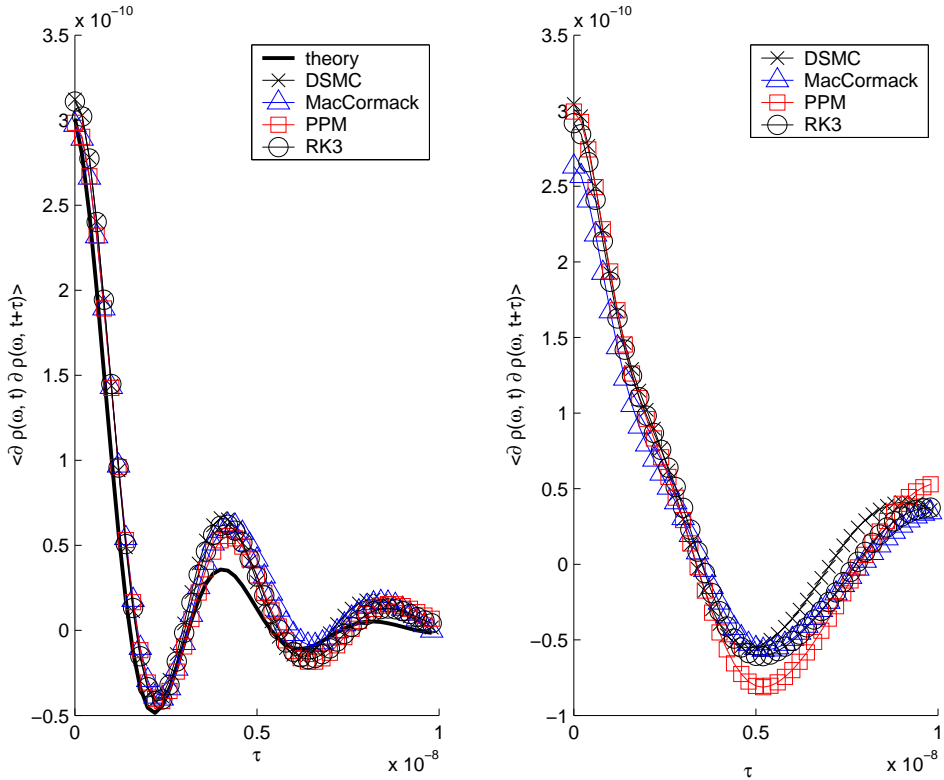}
\parbox{5in}{\caption{\label{TimeCorrPeriodicFig} Time correlation of density
fluctuations for equilibrium problem, on a periodic domain (left panel) and a
domain with specular wall boundaries (right panel).}}
\end{figure}

In the left-hand panel of figure~\ref{TimeCorrPeriodicFig}, we
present time correlation results from our equilibrium problem on a
periodic domain.  We compare results from the MacCormack, PPM, and
RK3 methods with the theoretical time correlation, equation
(\ref{TimeCorrelationEqn}), and with molecular simulation data (see
Appendix A). We find reasonable agreement among all the results, up
to the time when a sound wave has crossed the system ($\approx 4
\times 10^{-9}$ seconds). Due to finite size effects the theory is
only accurate for short times but the agreement among the numerical
PDE schemes and DSMC molecular simulation is good.

The right-hand panel of figure~\ref{TimeCorrPeriodicFig} shows time
correlation results for the equilibrium problem on a domain with
thermal walls rather than periodic boundaries; we find good
agreement for this problem as well, at least for times less than the
sound crossing time. For later times, the time correlation is
sensitive to the acoustic impedance of the thermal wall. For this
case, MacCormack under-predicts the correlation at early time while
PPM shows significant deviation near $t = 5 \times 10^{-8}$. Both
MacCormack and the RK3 scheme deviate somewhat from DSMC at late
time. Overall, however, the RK3 scheme captures the temporal
correlation better than either of the other two PDE schemes.

\section{Numerical Tests -- Non-equilibrium}\label{NumericalTestsNonEqSection}

The results from the section above indicate that of the three
stochastic PDE schemes, the third-order Runge-Kutta method (RK3)
consistently out-performs the other two schemes. In this section we
consider two more numerical tests, spatial correlations in a
temperature gradient and diffusion of a standing shock wave, but
restrict our attention to the RK3 scheme, comparing it with DSMC
molecular simulations.

\subsection{Spatial correlations in a temperature gradient}\label{SpaceCorrelationNonEqSection}

In the early 1980's, a variety of statistical mechanics calculations
predicted that a fluid under a non-equilibrium constraint, such as a
temperature gradient, would exhibit long-range correlations of
fluctuations~\cite{Schmitz:88}. Furthermore, quantities that are
independent at equilibrium, such as density and momentum
fluctuations, also have long-ranged correlations. These predictions
were qualitatively confirmed by light scattering
experiments~\cite{Beysens:80}, yet the effects are subtle and
difficult to measure accurately in the laboratory. Molecular
simulations confirm the predicted correlations of non-equilibrium
fluctuations for a fluid subjected to a temperature
gradient~\cite{Garcia:86,Mansour:87} and to a
shear~\cite{Garcia:87b}.

\begin{figure}
\centering\includegraphics[width=8.5cm]{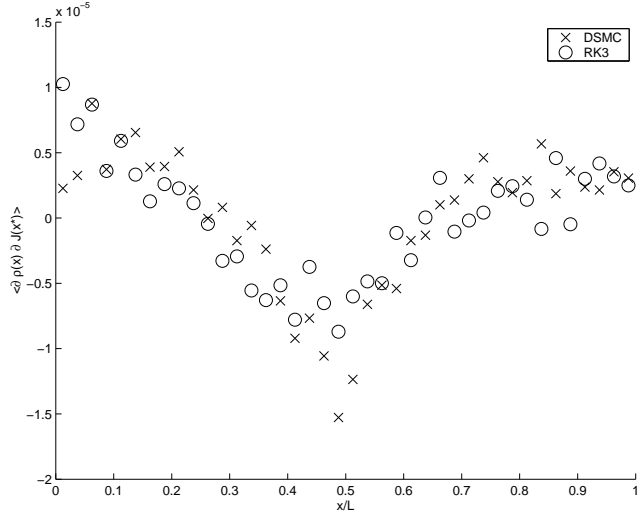}
\parbox{8.5cm}{\caption{\label{nonEqCorrRhoMomFig} Spatial correlation of density
and momentum fluctuations for a system subjected to a temperature
gradient. Compare with
Fig.~\ref{SpaceCorrRhoMomFig}.}}
\end{figure}

We consider a system similar to that of
section~\ref{TimeCorrelationSection}
but with a temperature
gradient. Specifically, the boundary conditions are thermal walls at
273K and 819K. Figure~\ref{nonEqCorrRhoMomFig} shows the correlation
of density and momentum fluctuations measured in an RK3 calculation
and by DSMC simulations. The two sets of data are in good agreement
and are in agreement with earlier work on this
problem~\cite{Garcia:86,Mansour:87}. The major discrepancy is the
under-prediction of the negative peak correlation near $j^*$.
Extensive tests suggest that this effect is hard to capture with a
continuum solver because of the tension between variance reduction
and spatial correlations in computing the mass flux at cell edges
from cell-centered data.

\subsection{Random Walk of a Standing
Shock}\label{ShockDiffusionSection}

In our final numerical study we consider the random walk of
a standing shock wave due to spontaneous fluctuations. Shock
diffusion is well-known in other particle simulations, such as shock
tube modeling by DSMC, which must correct for the drift when
measuring profiles for steady shocks.~\cite{Bird:94} The general problem has been
also been analyzed for simple lattice gas
models~\cite{Alexander:92,Alexander:93,Ferrari:94,Janowsky:92,Bell:06}.

Mass density and temperature on the right-hand side of the shock are
given the same values as in our equilibrium problem; values of
density and temperature on the left-hand side are derived from the
Rankine-Hugoniot relations. The velocity on both sides of the shock
are specified to satisfy the Rankine-Hugoniot conditions and to make
the unperturbed shock wave stationary in the computational domain.
We consider three different shock strengths, Mach 2, Mach 1.4, and
Mach 1.2 (see table \ref{standingShock_params}). The boundary
treatment consists of infinite reservoirs with the same states as
the initial conditions. For this test problem we use a longer
computational domain, in order to capture (unlikely) shock drift of
several standard deviations.

\begin{table}
  \centering
  \begin{tabular}{|l|c||l|c|}
    \hline
    System length & $5 \times 10^{-4}$      & Reference mean free path & $6.26  \times 10^{-6}$ \\
    System volume & $7.84 \times 10^{-16}$  & Time step & $1.0 \times 10^{-12}$\\
    Number of cells & 160                   & Mach number & 2.0 \\
\hline
    RHS mass density & $1.78\times 10^{-3}$ & LHS mass density & $4.07\times 10^{-3}$  \\
    RHS velocity     & -61562               & LHS velocity     & -26933 \\
    RHS temperature  & 273                  & LHS temperature  & 567    \\
    RHS sound speed  & 30781                & LHS sound speed  & 44373  \\
    \hline
  \end{tabular}
  \caption{System parameters (in cgs units) for simulations of a standing
  shock, Mach 2.0}\label{standingShock_params}
\end{table}

Here we focus on the variance of the shock location as a function of
time. We define a shock location for density, $\sigma_\rho(t)$ by
fitting a Heaviside function to the integrated density, i.e.,
\begin{equation}
\int_{-L/2}^{\sigma(t)} \rho_L \,dx + \int_{\sigma_\rho(t)}^{L/2}
\rho_R \,dx = \int_{-L/2}^{L/2} \rho(x,t) \,dx \;\;\; .
\end{equation}
Solving for $\sigma_\rho(t)$ gives
\begin{equation}
\sigma_\rho(t) = L\, \frac{\bar{\rho}(t) - \frac{1}{2} (\rho_L +
\rho_R)}{\rho_L - \rho_R}
\end{equation}
where $\bar{\rho} = L^{-1}\int_{-L/2}^{L/2} \rho(x,t) \,dx$ is the
instantaneous average density. The shock location for pressure,
$\sigma_P$, is analogously defined. We estimate $\sigma_\rho(t)$ and
$\sigma_p(t)$ as functions of time from ensembles of 4000
simulations. For the PDE simulations, we initialize with
discontinuous shock profiles. One would expect the shock location to
fluctuate with a diffusion similar to that of a simple random
walk~\cite{Ferrari:94}, so averaging over ensembles from the same
initial state we would expect to find
\begin{equation}
\langle \delta \sigma_\rho^2 \rangle \approx 2 \,\mathcal{D_\rho} t
 \qquad \mathrm{and} \qquad \langle \delta \sigma_p^2 \rangle \approx 2 \,\mathcal{D}_p t
\end{equation}
with shock diffusion coefficients, $\mathcal{D_\rho}$ and
$\mathcal{D}_p$, that depend on shock strength. Note that this
expression for the variance is not accurate at very short times (due
to transient relaxation from the initial state) or at very long
times (due to finite system size).

\begin{center}
\begin{figure}[h!]
\includegraphics[width=3.0in]{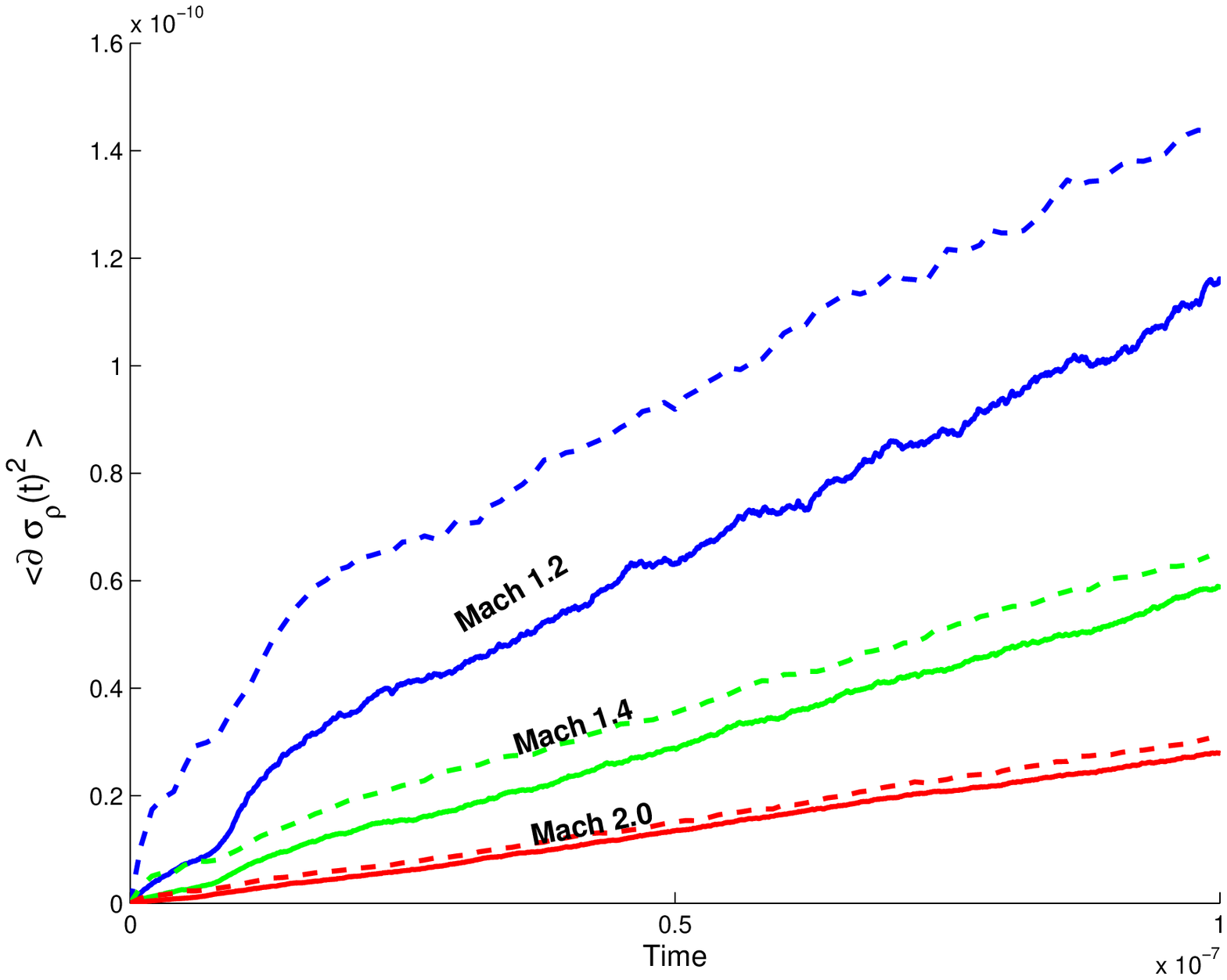}\includegraphics[width=3.0in]{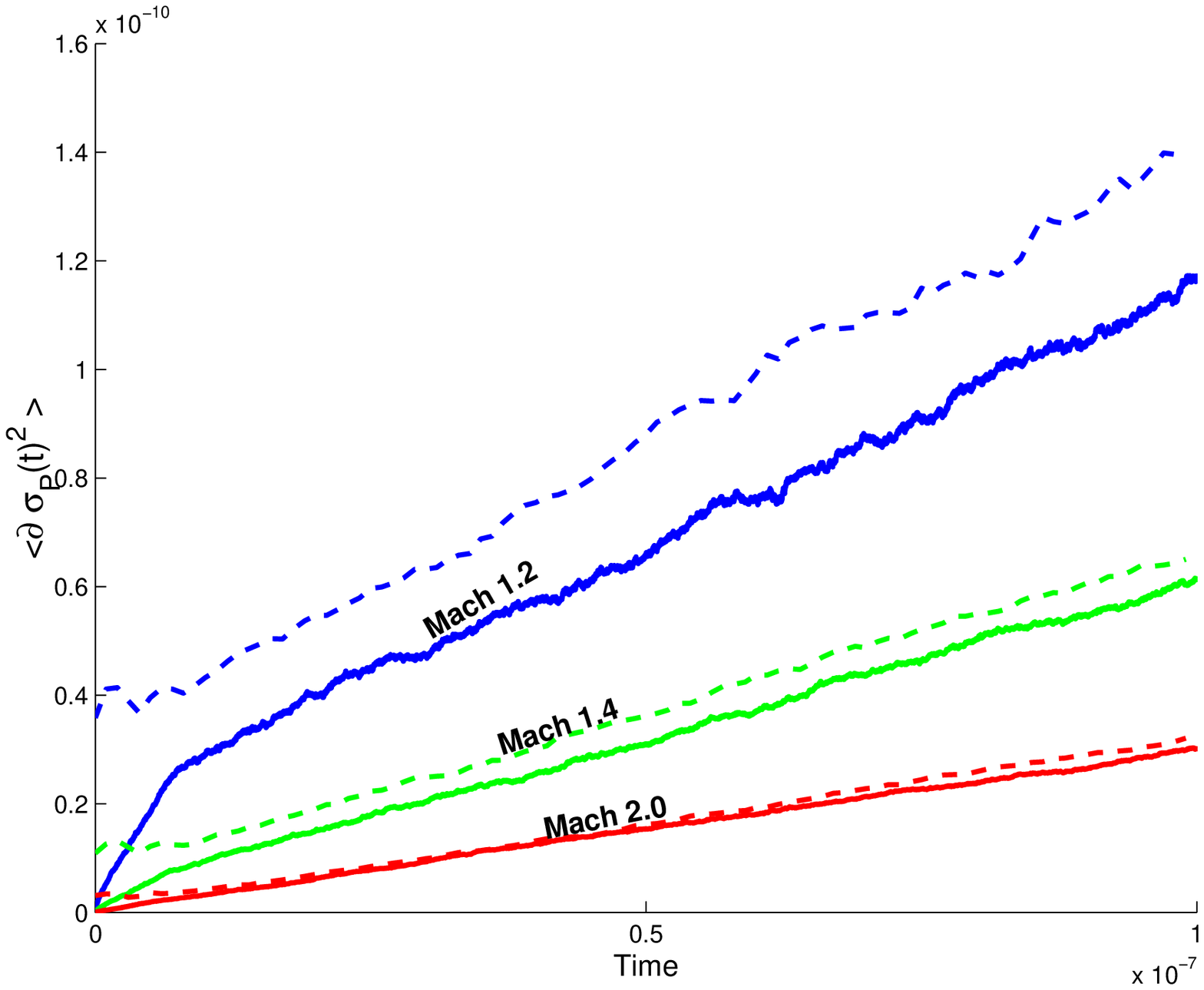}
\caption{Variance of shock location for mass density profile (left panel, $\langle \delta
  \sigma_{\rho}(t)^2 \rangle$)
and pressure profile (right panel, $\langle \delta
  \sigma_{P}(t)^2 \rangle$).  Estimated variances
(4000-run ensembles) versus time $t$ for a deterministically steady shock of Mach number
1.2, 1.4, or 2.0. Solid lines are for RK3, dashed lines
are from DSMC molecular simulations.}\label{ShockDiffusionFig}
\end{figure}
\end{center}

Figure~\ref{ShockDiffusionFig} shows results for the variance in the
shock position from an ensemble of runs versus time. After the
initial transients, the slopes are constant with the strongest
shocks exhibiting the least drift ($\mathcal{D} \cong
(\mathrm{Ma}-1)^{-1}$) and with $\sigma_\rho$ and $\sigma_P$ giving
similar diffusion coefficients. DSMC data is initially noisy so it
has different initial transients and ``diffuses'' farther than the
PDE. However, after the transients, the DSMC and the RK3 simulations
have essentially the same slope, as a function of Mach number. This
indicates that the third-order Runge-Kutta scheme is accurately
capturing the shock-drift random walk.

\section{Summary and Concluding Remarks}\label{SummarySection}

In this paper we develop and analyze several finite-volume schemes
for solving the fluctuating Landau-Lifshitz compressible
Navier-Stokes equations in one spatial dimension. Methods based on
standard CFD discretizations were found not to accurately represent
fluctuations in an equilibrium flow. We have introduced a centered
scheme based on interpolation schemes designed to preserve
fluctuations combined with a third-order Runge-Kutta (RK3) temporal
integrator that was able to capture the equilibrium fluctuations.
Further tests for non-equilibrium systems confirm that the RK3 scheme
correctly reproduces long-ranged correlations of fluctuations and
stochastic drift of shock waves, as verified by comparison with
molecular simulations. It is worth emphasizing that the ability of
continuum methods to accurately capture fluctuations is fairly
sensitive to the construction of the numerical scheme. Minor
variations in the numerics can lead to significant changes in
stability, accuracy, and behavior.

The work discussed here suggests a number of additional studies.
Further analysis is needed on the treatment of thermal and reservoir
boundary conditions. The methods here can also be extended to three
dimensions (for which the stochastic stress tensor is more complex)
and we can include concentration as a hydrodynamic variable to allow
the methodology to be applied to a number of other flow problems.
Finally, we are integrating our new stochastic PDE solver into our
existing Adaptive Mesh and Algorithm Refinement (AMAR)
programs~\cite{Garcia:99}. A stochastic AMAR simulation will not
only model hydrodynamic fluctuations at multiple grid scales but
will, by incorporating DSMC simulations at the finest level of
algorithm refinement, also capture molecular-level physics.

\section* {Acknowledgment}

The authors wish to thank Phil Colella for
helpful discussions about the PPM.
The work of John Bell was supported by
the Applied Mathematics Program of the DOE Office of Mathematics,
Information, and Computational Sciences under the U.S. Department of
Energy under contract No.\ DE-AC03-76SF00098.
Sarah Williams' support was provided by
DE-FC02-01ER25473 SciDAC and DE-FG02-03ER25579 MICS grants.

%
%

\section*{Appendix A: Equilibrium Fluctuations}

For infinite systems, at thermodynamic equilibrium both conserved and hydrodynamic
variables are spatially uncorrelated at equal times.  For example,
\begin{equation}
\langle \delta \rho_i(t) \delta \rho_j(t) \rangle = \langle \delta
\rho^2 \rangle \delta^K_{i,j}.
\label{EquilibCorrRhoEqn}
\end{equation}
For conserved variables there is a finite size correction, specifically,
\begin{equation}
\langle \delta \rho_i(t) \delta \rho_j(t) \rangle =
(1 - M_c^{-1}) \langle \delta \rho^2 \rangle \delta^K_{i,j}
 - M_c^{-1} \langle \delta \rho^2 \rangle (1 - \delta^K_{i,j})
\label{EquilibCorrRhoFiniteEqn}
\end{equation}
for $i,j = 1,\ldots,M_c$, where $M_c$ is the number of cells in the system.
The variances are well-known from equilibrium statistical mechanics
(\S 112, \cite{Landau:StatPhys1}).

The variance of mass density depends on the compressibility (i.e., the equation of state) of the
fluid. In general,
\begin{equation}
\langle{ \delta \rho^2 }\rangle =
\overline{\rho}^2~\frac{\langle \delta N_c^2 \rangle}{\overline{N}_c^2}
\label{EquilibVarRhoEqn}
\end{equation}
where $\overline{N}_c$ and $\langle \delta N_c^2 \rangle$ are the mean and
variance of the number of particles in a cell.
We calculate $\overline{N}_c = \overline{\rho} V_c/m$, where $V_c$ is the volume
of a cell and $m$ is the mass of a particle.
For an ideal gas $N_c$ is Poisson
distributed so $\langle \delta N_c^2 \rangle = \overline{N}_c$ and $\langle \delta
\rho^2 \rangle = \overline{\rho}^2/\overline{N}_c$. The more general result
is $\langle \delta N_c^2 \rangle = \alpha_T \rho k_B \overline{T}\,\overline{N}_c/m$ where
$\alpha_T$ is the isothermal compressibility.

The variances of fluid velocity and temperature in a cell are
\begin{eqnarray}
\langle{ \delta u^2 }\rangle &=& \frac{k_B
\overline{T}}{\overline{\rho} V_c}
= \frac{C_T^2}{\overline{N}_c}\label{EquilibVarVelEqn} \textrm{ and}\\
\langle{ \delta T^2 }\rangle &=& \frac{k_B \overline{T}^2}{c_v
\overline{\rho} V_c} = \frac{C_T^2 \overline{T}}{c_v
\overline{N}_c}\label{EquilibVarTempEqn},
\end{eqnarray}
where $C_T = \sqrt{k_B \overline{T}/m}$ is the thermal speed (and
the standard deviation of the Maxwell-Boltzmann distribution).
The
covariances are $\langle{ \delta \rho \,\delta u }\rangle =
\langle{ \delta \rho \,\delta T }\rangle =
\langle{\delta u\,\delta T}\rangle = 0$.

The variances and covariances of the mechanical densities at
equilibrium are
\begin{eqnarray}
\langle{ \delta \rho \delta J }\rangle &=&
\overline{\rho} \overline{J} \Delta_\rho  \label{EquilibCorrRhoMomEqn}\\
\langle{ \delta \rho \delta E }\rangle &=&
\overline{\rho} \overline{E} \Delta_\rho  \label{EquilibCorrRhoEnrEqn}\\
\langle{ \delta J^2}\rangle &=&
\overline{J}^2  \Delta_\rho
+ \overline{\rho}^2 C_T^2 \Delta_u   \label{EquilibVarMomEqn}\\
\langle{ \delta J\, \delta E }\rangle &=&
\overline{J}\,\overline{E} \Delta_\rho
+ \overline{J}\,\rho C_T^2 \Delta_u  \label{EquilibCorrMomEnrEqn}\\
\langle{ \delta E^2 }\rangle &=& \overline{E}^2 \Delta_\rho +
\overline{J}^2 C_T^2 \Delta_u + c_v^2 \overline{\rho}^2
\overline{T}^2 \Delta_T  \label{EquilibVarEnrEqn}
\end{eqnarray}
where $\Delta_\rho = {\langle \delta \rho^2
\rangle}/{\overline{\rho}^2}$, $\Delta_u = {\langle \delta u^2
\rangle}/{C_T^2}$, and $\Delta_T = {\langle \delta T^2
\rangle}/{\overline{T}^2}$.
For a dilute gas $\Delta_\rho = \Delta_u = 1/\overline{N}_c$, and $\Delta_T = 2/(3\overline{N}_c)$.
Again,
corrections must be made for conserved quantities in the case of a finite domain:
\begin{equation}
\langle \delta J_i(t) \delta J_j(t) \rangle =
(1 - M_c^{-1}) \langle \delta J^2 \rangle \delta^K_{i,j}
 - M_c^{-1} \langle \delta J^2 \rangle (1 - \delta^K_{i,j}),
\label{EquilibCorrMomFiniteEqn}
\end{equation}
\begin{equation}
\langle \delta E_i(t) \delta E_j(t) \rangle =
(1 - M_c^{-1}) \langle \delta E^2 \rangle \delta^K_{i,j}
 - M_c^{-1} \langle \delta E^2 \rangle (1 - \delta^K_{i,j}).
\label{EquilibCorrEnrFiniteEqn}
\end{equation}
%

\section*{Appendix B: DSMC Simulations}

The algorithms presented here for the stochastic LLNS equations were
validated by comparison with molecular simulations. Specifically, we
used the direct simulation Monte Carlo (DSMC) algorithm, a
well-known method for computing gas dynamics at the molecular scale;
see \cite{Alexander:97,Garcia:00} for pedagogical expositions on
DSMC, \cite{Bird:94} for a complete reference, and \cite{Wagner:92}
for a proof of the method's equivalence to the Boltzmann equation.
As in molecular dynamics, the state of the system in DSMC is given
by the positions and velocities of particles. In each time step, the
particles are first moved as if they did not interact with each
other. After moving the particles and imposing any boundary
conditions, collisions are evaluated by a stochastic process,
conserving momentum and energy and selecting the post-collision
angles from their kinetic theory distributions. DSMC is a stochastic
algorithm but the statistical variation of the physical quantities
has nothing to do with the ``Monte Carlo'' portion of the method.
The equilibrium and non-equilibrium variations in DSMC are the
physical spectra of spontaneous thermal fluctuations, as confirmed
by excellent agreement with fluctuating hydrodynamic
theory~\cite{Garcia:87,Mansour:87} and molecular dynamics
simulations~\cite{Mansour:88,Mareschal:92}.

The simulated physical system is a dilute monatomic hard-sphere gas
in a rectangular volume with periodic boundary conditions in the $y$
and $z$ directions. The boundary conditions in the $x$ direction are
either periodic, specular (i.e., elastic reflection of particles),
or a pair of parallel thermal walls. The physical parameters used
are presented in Table~\ref{DSMC_Table}. Samples are taken in forty
rectangular cells perpendicular to the $x$-direction.

\bibliographystyle{unsrt}

\bibliography{fluct_hydro_bib}

\end{document}